# Bootstrapping the Grenander estimator

## Michael R. Kosorok[*1]

*University of North Carolina-Chapel Hill*

**Abstract:** The goal of this paper is to study the bootstrap for the Grenander
estimator. The first result is a proof of the inconsistency of the nonparametric
bootstrap for the Grenander estimator at a given point. The second result is
the development and verification of a bootstrap for the $L_1$ confidence band
for the Grenander estimator. As part of this work, kernel estimators are studied
as alternatives to the Grenander estimator. We show that when the second
derivative of the true density is assumed to be uniformly bounded, there exist
kernel estimators with faster convergence rates than the Grenander estimator.
We study the implications of this in developing $L_1$ and uniform confidence
bands and discuss some open questions.

## 1. Introduction

The Grenander estimator (Grenander [6]) is the maximum likelihood estimator
(MLE) $\hat{f}_n$ of the density $f$ of a positive, real random variable $X$, under the con-
straint that $f$ is monotone non-increasing. For simplicity, we will assume the support
of $X$ is $[0, 1]$, and that the data used for estimation is an i.i.d. sample $X_1, \ldots, X_n$
from $f$.

Not only is the Grenander estimator worthy of study in its own right, but it
is also useful because of its connection to the MLE of the survival function from
current status time-to-event data. Current status data arises from only observing
the current status of the event time $T$ at a random observation time $Y$. Specifically,
one only observes $\mathbf{1}\{T \leq Y\}$ and $Y$, where $\mathbf{1}\{A\}$ denotes the indicator of $A$. The
connection between the estimators is that the current status MLE of the cumulative
distribution at a single point and $\hat{f}_n(t)$ for a point $t$ both have the same limiting
distribution after suitable standardization. This limiting distribution is

$$\mathbb{C} \equiv \underset{h \in \mathbb{R}}{\arg\max} \left\{ \mathbb{Z}(h) - h^2 \right\},$$

where $\mathbb{Z}$ is two-sided Brownian motion with $\mathbb{Z}(0) = 0$. This result was obtained for
the Grenander estimator by Rao [16] (see also Groeneboom [7]) and for the current
status estimator by Groeneboom [8] (see also Groeneboom and Wellner [10]).

Because of this similarity between the Grenander estimator and the current
status survival estimator, it appears that at least some of what we could learn about
the Grenander estimator may be applicable on many levels to inference problems

---

[*]Supported in part by Grant CA075142 from the National Cancer Institute.
[1]Department of Biostatistics, School of Public Health, University of North Carolina-Chapel
Hill, 3101 McGavran-Greenberg Hall, CB 7420, Chapel Hill, NC 27599-7420, USA, e-mail:
kosorok@unc.edu

in current status and possibly other, more complex, data types. We will not pursue this connection further in this paper, but the interested reader should compare and contrast the derivations of both of these estimators as given side-by-side in Sections 3.2.14 and 3.2.15 of van der Vaart and Wellner [18].

There are number of challenges with using $\mathbb{C}$ – Chernoff's distribution – directly for inference. The first challenge is that the density for $\mathbb{C}$, the form of which was derived by Groeneboom [9], does not have a closed form. Computing critical values for the distribution is quite difficult (Dykstra and Carolan [5], Narayanan and Sager [14]) but has been done (Groeneboom and Wellner [11]). A second challenge is that the normalization constants involved can be difficult to estimate. For these reasons, computationally reasonable approaches that avoid computing critical values from Chernoff's distribution, and/or ameliorate the need for computing complicated normalization constants, would be very appealing. This pursuit is the theme of this paper.

An obvious approach to consider because of computational simplicity is the nonparametric bootstrap. Unfortunately, the first main result of this paper is that the nonparametric bootstrap is inconsistent for pointwise inference (i.e., inference for $f(t)$ at a given value of $t \in [0, 1]$). We prove this rigorously in Theorem 2.1 below. The key argument is contained in Theorem 2.2 below and is applicable to many other inference settings. The inconsistency of the bootstrap was also observed by Abrevaya and Huang [1] for the maximum score estimator, which also has a Chernoff limit, although they did not provide a rigorous proof of this. Fortunately, we are able to show that a smoothed bootstrap (Silverman and Young [17]) obtained by sampling from a certain kernel estimator is consistent. These pointwise inference results will be presented in Section 2.

It would be nice if some of the pointwise results could be utilized in the development of uniform confidence bands, but this appears to be an excruciatingly difficult problem. However, some progress has been made for $L_1$ confidence bands. Building on the work of Groeneboom et al. [12], who derive the limiting distribution of the $L_1$ error of the Grenander estimator, we propose a "supersampling" smoothed bootstrap. This is discussed in Section 3. One of the discoveries made in this process is that the assumptions needed for $L_1$ convergence of the Grenander estimator are so strong that there exist kernel estimators with faster convergence rates than the Grenander estimator.

We conclude the paper with a discussion of the implication of these results and several open questions in Section 4. The main contributions of the paper are first, a proof of the invalidity of the nonparametric bootstrap for the Grenander estimator, and, second, the development of smoothed bootstrap procedures for both pointwise and $L_1$ confidence bands. The results and ideas of this paper should prove useful in developing solutions to the confidence band problem for the Grenander estimator as well as for current status survival function estimators and other related monotone function estimators.

## 2. Pointwise error

The focus of this section is on pointwise inference based on the Grenander estimator. Before presenting the main results on the bootstrap, we first briefly review known asymptotic distribution results for the Grenander estimator $\hat{f}_n$. Before doing this, however, we make the following assumptions about the density $f$:

A1. $0 \leq f(1) \leq f(s) \leq f(t) \leq f(0) < \infty$, for all $0 \leq t \leq s \leq 1$; and

A2. $f$ is differentiable with derivative $\dot{f}$ satisfying

$$0 < \inf_{t \in (0,1))} |\dot{f}(t)| \leq \sup_{t \in (0,1)} |\dot{f}(t)| < \infty.$$

We may occasionally need stronger assumptions which will be introduced as needed.

It is well known that $\hat{f}_n$ is the left derivative of the least concave majorant of the empirical distribution function $\mathbb{F}_n$ (see, for example, Section 3.2.14 of van der Vaart and Wellner [18]). Moreover, under assumptions A1 and A2, we have the now classic result that for any $t \in (0,1)$,

$$n^{1/3}(\hat{f}_n(t) - f(t)) \quad \rightsquigarrow \quad |4\dot{f}(t)f(t)|^{1/3} \arg\max_{h \in \mathbb{R}} \left\{ \mathbb{Z}(h) - h^2 \right\}$$
$$\equiv \quad c(t)\mathbb{C}$$

(Groeneboom [8]). Both the normalizing constant $c(t)$ and critical values for Chernoff's distribution are needed for inference.

Before proceeding to the bootstrap discussion, we point out that an alternative to the above approach to inference about $f$ is to use the nonparametric maximum likelihood ratio as done by Banerjee and Wellner [2] which has an asymptotically pivotal distribution that avoids the need to estimate a normalizing constant. The limiting distribution in this setting is not Chernoff's distribution but is still quite complicated and does not have a known, closed form. Computing critical values is possible but complicated (Banerjee and Wellner [3]).

Another alternative that almost always works theoretically is the subsampling bootstrap of Politis and Romano [15]. The basic idea is to perform a bootstrap without replacement of sample size $m$ which is much smaller than the actual sample size $n$. Provided $m \to \infty$ and $m/n \to 0$, the standardized subsample bootstrap will be valid. Unfortunately, in practice, this will not work unless $n$ is quite large since the asymptotic approximation must be approximately valid for the subsample size $m$, not just valid for $n$. We will not pursue the subsample bootstrap further in this paper.

Let $\mathbb{F}_n^*$ be the usual nonparametric bootstrap empirical distribution function, and let $\hat{f}_n^*$ be the the left derivative of the least concave majorant of $\mathbb{F}_n^*$. What we would like to show is that $n^{1/3}(\hat{f}_n^*(t) - \hat{f}_n(t))$, conditional on the data $X_1, X_2, \ldots,$ converges to the unconditional limiting distribution of $n^{1/3}(\hat{f}(t) - f(t))$. Our first main result is that this approach is unfruitful, as we now show in the following theorem:

**Theorem 2.1.** *The nonparametric bootstrap is inconsistent for the Grenander estimator, i.e., $n^{1/3}(\hat{f}_n^*(t) - \hat{f}_n(t))$ does not converge in probability, conditional on the data, to $c(t)\mathbb{C}$, for any $t \in (0,1)$.*

Before giving the proof of this theorem, we present a general theorem which can be useful in studying bootstrap validity. Let $X_n$ be a random variable in a Banach space $(\mathbb{B}, \|\cdot\|)$ that converges weakly to a tight limit $X$, and let $\hat{X}_n$ denote a bootstrapped version of $X_n$ based on some random weighting mechanism $W_n$ which is independent of the data $\mathcal{X}_n$ used to generate $X_n$. We say that $\hat{X}_n$ is a valid bootstrap if its limiting distribution conditional on $\mathcal{X}_n$ "converges weakly" to $X$.

We now define what "converges weakly" means in this context. Let $BL_1(\mathbb{B})$ be the collection of all Lipschitz continuous functions $h : \mathbb{B} \mapsto \mathbb{R}$ bounded in absolute value by 1 and having Lipschitz constant 1, i.e., $|h| \leq 1$ and $|h(x) - h(y)| \leq \|x - y\|$ for all $x, y \in \mathbb{B}$. We say $\hat{X}_n$ converges weakly conditional on the data to $X$ if

$$\sup_{h \in BL_1(\mathbb{B})} \left| \mathrm{E}_{\cdot|\mathcal{X}_n} h(\hat{X}_n) - \mathrm{E}h(X) \right| \to 0,$$

in outer probability, where $\mathrm{E}_{\cdot|\mathcal{X}_n}$ denotes conditional expectation given $\mathcal{X}_n$, and provided $h(\hat{X}_n)$ is asymptotically measurable unconditionally for all $h \in BL_1(\mathbb{B})$. We denote this kind of conditional convergence $\hat{X}_n \underset{W}{\overset{\mathrm{P}}{\rightsquigarrow}} X$. We also require $h(\hat{X}_n)$ to be a measurable function of $W_n$ conditional on $\mathcal{X}_n$ for all $h \in BL_1(\mathbb{B})$. A more precise discussion of this general formulation of the bootstrap can be found in van der Vaart and Wellner [18]. The following is a general result for these kinds of bootstraps:

**Theorem 2.2.** *Assume $X_n \rightsquigarrow X$, where $X$ is tight, and that $\hat{X}_n \underset{W}{\overset{\mathrm{P}}{\rightsquigarrow}} X$, where $W_n \mapsto h(\hat{X}_n)$ is measurable conditional on $\mathcal{X}_n$ for all $h \in BL_1(\mathbb{B})$. Then $(\hat{X}_n, X_n) \rightsquigarrow (\tilde{X}_1, \tilde{X}_2)$ unconditionally, where $\tilde{X}_1$ and $\tilde{X}_2$ are independent copies of $X$.*

*Proof.* Let $\tilde{X}_1$ and $\tilde{X}_2$ be two independent copies of $X$, which are also independent of the data $\mathcal{X}_n$, and note that

$$\sup_{h \in BL_1(\mathbb{B}^2)} \left| \mathrm{E}^* h(\hat{X}_n, X_n) - \mathrm{E}h(\tilde{X}_1, \tilde{X}_2) \right|$$
$$\leq \sup_{h \in BL_1(\mathbb{B}^2)} \left| \mathrm{E}^* h(\hat{X}_n, X_n) - \mathrm{E}^* h(\tilde{X}_1, X_n) \right|$$
$$+ \sup_{h \in BL_1(\mathbb{B}^2)} \left| \mathrm{E}^* h(\tilde{X}_1, X_n) - \mathrm{E}h(\tilde{X}_1, \tilde{X}_2) \right|$$
$$\equiv A_n + B_n.$$

Note also that for any $h \in BL_1(\mathbb{B}^2)$ and any $y \in \mathbb{B}$, both $x \mapsto h(x, y)$ and $x \mapsto h(y, x)$ are members of $BL_1(\mathbb{B})$. As a consequence of the weak convergence of $\hat{X}_n$, we therefore have

$$\sup_{h \in BL_1(\mathbb{B}^2)} \left| \mathrm{E}_{\cdot|\mathcal{X}_n} h(\hat{X}_n, X_n) - \mathrm{E}_{\cdot|\mathcal{X}_n} h(\tilde{X}_1, X_n) \right|$$
$$\leq \sup_{h \in BL_1(\mathbb{B})} \left| \mathrm{E}_{\cdot|\mathcal{X}_n} h(\hat{X}_n) - \mathrm{E}h(\tilde{X}_1) \right|$$
$$\overset{\mathrm{P}}{\to} 0,$$

where $\overset{\mathrm{P}}{\to}$ denotes convergence in outer probability. Provided both

$$(2.1) \qquad \sup_{h \in BL_1(\mathbb{B}^2)} \left| \mathrm{E}^* \mathrm{E}_{\cdot|\mathcal{X}_n} h(\hat{X}_n, X_n) - \mathrm{E}^* h(\hat{X}_n, X_n) \right| \quad \to \quad 0$$

and

$$(2.2) \qquad \sup_{h \in BL_1(\mathbb{B}^2)} \left| \mathrm{E}^* \mathrm{E}_{\cdot|\mathcal{X}_n} h(\tilde{X}_1, X_n) - \mathrm{E}^* h(\tilde{X}_1, X_n) \right| \quad \to \quad 0,$$

we will obtain that $A_n \to 0$.

Arguing in a similar manner but utilizing instead the assumed weak convergence of $X_n$, we have

$$\sup_{h \in BL_1(\mathbb{B}^2)} \left| \mathrm{E}^*_{\cdot|\tilde{X}_1} h(\tilde{X}_1, X_n) - \mathrm{E}_{\cdot|\tilde{X}_1} h(\tilde{X}_1, \tilde{X}_2) \right| \quad \leq \quad \sup_{h \in BL_1(\mathbb{B})} \left| \mathrm{E}^* h(X_n) - \mathrm{E}h(\tilde{X}_2) \right|$$
$$\to \quad 0,$$

where $\mathrm{E}_{\cdot|\tilde{X}_1}$ denotes conditional expectation given $\tilde{X}_1$. Provided

$$(2.3) \qquad \sup_{h \in BL_1(\mathbb{B}^2)} \left| \mathrm{E}^* \mathrm{E}_{\cdot|\tilde{X}_1}^* h(\tilde{X}_1, X_n) - \mathrm{E}^* h(\tilde{X}_1, X_n) \right| \to 0,$$

we will obtain that $B_n \to 0$, and the desired conclusion of the theorem will follow.

The proof is essentially complete, except for establishing (2.1), (2.2), and (2.3), which are primarily measurability technicalities. The uninterested reader can skip this part of the proof and proceed directly to the proof of Theorem 2.1 below. Since $h(\hat{X}_n)$ is asymptotically measurable for all $BL_1(\mathbb{B})$, it is asymptotically measurable for all $h$ that is bounded and Lipschitz continuous for any Lipschitz constant. Thus the conditional weak convergence of $\hat{X}_n$ implies that for every $\epsilon > 0$, there exists a compact $K \subset \mathbb{B}$ such that

$$\mathrm{E}^* \mathrm{P} \left( \hat{X}_n \in K^\delta \,\middle|\, \mathcal{X}_n \right) \to P(X \in K^\delta) \geq 1 - \epsilon,$$

for every $\delta > 0$, where $K^\delta \equiv \{x \in \mathbb{B} : \|x - y\| < \delta, \text{ for some } y \in K\}$. Hence, by Fubini's theorem for outer expectations (Lemma 1.2.6 of van der Vaart and Wellner [18]), $\hat{X}_n$ is asymptotically tight unconditionally. Thus it is also asymptotically measurable unconditionally by reapplication of Lemma 1.3.13 of van der Vaart and Wellner [18]. Since marginal asymptotic tightness plus marginal asymptotic measurability implies joint asymptotic tightness and measurability (see Lemmas 1.4.3 and 1.4.4 of van der Vaart and Wellner [18]), we have that $(\hat{X}_n, X_n)$ is jointly asymptotically tight and measurable. Thus

$$\sup_{h \in BL_1(\mathbb{B}^2)} \left| \mathrm{E}^* h(\hat{X}_n, X_n) - \mathrm{E}_* h(\hat{X}_n, X_n) \right| \to 0,$$

and condition (2.1) follows.

The assumed weak convergence of $X_n$ implies asymptotic measurability via Lemma 1.3.13 of van der Vaart and Wellner [18], and thus (2.2) also follows. Since $X_n$ converges weakly, $(\tilde{X}_1, X_n)$ jointly converges weakly, and thus $(\tilde{X}_1, X_n)$ is asymptotically measurable. Hence

$$\sup_{h \in BL_1(\mathbb{B}^2)} \left| \mathrm{E}^* h(\tilde{X}_1, X_n) - \mathrm{E}_* h(\tilde{X}_1, X_n) \right| \to 0,$$

and (2.3) will follow. This completes the proof in all of its formality. $\qquad \square$

*Proof of Theorem 2.1.* The basic idea of the proof is to assume that

$$(2.4) \qquad n^{1/3}(\hat{f}_n^*(t) - \hat{f}_n(t)) \quad \underset{W}{\overset{\mathrm{P}}{\rightsquigarrow}} \quad c(t)\mathbb{C},$$

where the $W$ refers to the random multinomial weights $W_n \equiv \{W_{n,1}, \ldots, W_{n,n}\}$ in the nonparametric bootstrap, and then use Theorem 2.2 to obtain a contradiction. Accordingly, assume (2.4), and let $\hat{X}_n = n^{1/3}(\hat{f}_n^*(t) - \hat{f}_n(t))$ and $X_n = n^{1/3}(\hat{f}_n(t) - f(t))$. Then Theorem 2.2 implies that $\hat{X}_n + X_n \rightsquigarrow c(t)(\mathbb{C}_1 + \mathbb{C}_2)$, unconditionally, where $\mathbb{C}_1$ and $\mathbb{C}_2$ are two independent copies of $\mathbb{C}$.

Since $Y_n \equiv \hat{X}_n + X_n = n^{1/3}(\hat{f}_n^* - f(t))$, the above results imply that $Y_n$ converges unconditionally to a tight limiting distribution which has twice the variance of $c(t)\mathbb{C}$. Using arguments along the lines of those used in Section 3.2.14 of van der Vaart and Wellner [18], along with properties of bootstrapped empirical processes, it is not hard to verify, however, that

$$
\left[\frac{n}{4\dot{f}(t)f(t)}\right]^{1/3}(\hat{f}_n^*(t) - f(t)) \quad \rightsquigarrow \quad \operatorname{argmax}_h\left\{\mathbb{Z}_1(h) + \mathbb{Z}_2(h) - h^2\right\}
$$

(2.5)
$$
\equiv \quad \tilde{\mathbb{C}},
$$

where $\mathbb{Z}_1$ and $\mathbb{Z}_2$ are independent two-sided Brownian motions.

Using symmetry properties of Brownian motion and a careful change of variables, we can derive that $\tilde{\mathbb{C}}$ has the same distribution as $a \arg\max_h\{\sqrt{2a}\mathbb{Z}(ah) - (ah)^2\}$, for a two-sided Brownian motion $\mathbb{Z}$. Choosing $a = 2^{1/3}$ yields that $\tilde{\mathbb{C}}$ has the same distribution as $2^{1/3}\mathbb{C}$, and thus the variance of the limiting distribution of $Y_n$ is $2^{2/3} < 2$ times the variance of $c(t)\mathbb{C}$. This is a contradiction, and thus the desired conclusion follows. □

We now work toward developing an asymptotically valid alternative to the nonparametric bootstrap. To accomplish this, we propose a version of the "smoothed" bootstrap (Silverman and Young [17]). The idea is that we estimate the density with a certain modified kernel density estimator $\check{f}_n$, and then draw a smoothed bootstrap sample from $\tilde{f}_n$. Our goal is to ensure that the properties of this procedure lead to valid inference.

Let the kernel be $K$ and assume the bandwidth $1/2 \geq h \to 0$ as $n \to \infty$. For all $t \in [h, 1-h]$, let

$$
\check{f}_n = \int_0^1 \frac{1}{h}K\left(\frac{t-u}{h}\right)d\mathbb{F}_n(u),
$$

and denote $\check{f}_n^{(1)}$ as the first derivative of $\check{f}_n$ (so far only defined on $[h, 1-h]$). For $t \in [0, h)$, let

$$
\check{f}_n(t) = \check{f}_n(h) + (t-h)\left\{\check{f}_n^{(1)}(h) \wedge 0\right\},
$$

and for $t \in (1-h, 1]$, let

$$
\check{f}_n(t) = \check{f}_n(1-h) + (t-1+h)\left\{\check{f}_n^{(1)}(1-h) \wedge 0\right\}.
$$

Finally, define

$$
\tilde{f}_n(t) = \frac{\check{f}_n(t) \vee 0}{\int_0^1\left\{\check{f}_n(s) \vee 0\right\}ds}.
$$

We need the following assumptions on $K$:

B1. The kernel $K$ is nonnegative with support on $[-1, 1]$;

B2. $K$ is bounded and $\int_{-1}^1 K(v)dv = 1$;

B3. $\dot{K}$ is bounded, $v\dot{K}(v) \leq 0$ for all $v \in [-1, 1]$, $\int_{-1}^1 \dot{K}(v)dv = 0$, and $\int_{-1}^1 v\dot{K}(v)dv = -1$; and

B4. $|\ddot{K}|$ is uniformly bounded over $(-1, 1)$.

Two examples of kernels that satisfy B1–B4 are $K(v) = (3/4)(1-v^2)$ and $K(v) = (15/16)(1-v^2)^2$.

We now have the following lemma:

**Lemma 2.1.** *Provided $h = R_n n^{-\alpha}$, where $0 < R_n + R_n^{-1} = O_P(1)$ and $\alpha \in (0, 1/3)$, we have the following under assumptions A1–A2 and B1–B4:*

(i) $\tilde{f}_n$ *is uniformly consistent for* $f$;

(ii) *There exists constants* $0 < a < b < \infty$ *such that*

$$-b - o_P(1) \leq \inf_{t \in (0,1)} \tilde{f}_n^{(1)}(t) \leq \sup_{t \in (0,1)} \tilde{f}_n^{(1)}(t) \leq -a + o_P(1).$$

*Under the additional assumption*

A3. $\dot{f}(t)$ *is continuous at* $t = t_0$, *for some* $t_0 \in (0,1)$,

*we also have*

(iii) $\tilde{f}_n^{(1)}(t_0) = \dot{f}(t_0) + o_P(1)$.

*Proof.* The proof of conclusion (i) follows from standard arguments, and we omit the details. For (ii), we use change of variables to obtain that

$$
\begin{aligned}
\int_{-1}^{1} h^{-2} \dot{K}\left(\frac{t-u}{h}\right) d\mathbb{F}_n(u) - \dot{f}(t) &= n^{-1/2} \int_{-1}^{1} h^{-2} \dot{K}\left(\frac{t-u}{h}\right) d\mathbb{H}_n(u) \\
&\quad + h^{-1} \int_{-1}^{1} \dot{K}(v)(f(t-hv) - f(t)) dv \\
&= O_P\left(n^{-1/2} h^{-2} \sup_{n|s-t| \leq h} |\mathbb{G}_n(s) - \mathbb{G}_n(t)|\right) \\
&\quad + h^{-1} \int_{0}^{1} \dot{K}(v)\left[-\dot{f}(t_{hv}) h v\right] dv,
\end{aligned}
$$

where $u \mapsto \mathbb{H}_n(u) \equiv \mathbb{G}_n(u) - \mathbb{G}_n(t)$, $\mathbb{G}_n \equiv \sqrt{n}(\mathbb{F}_n - F)$ and $t_{hv}$ is on the line segment between $t$ and $t - hv$. By A1 and A2 combined with the fact that

$$n^{-1/2} h^{-3/2} \sqrt{\log\left(\frac{1}{h}\right)} = o(1),$$

conclusion (ii) follows. Conclusion (iii) follows because, when $\dot{f}$ is continuous at $t_0$, $\dot{f}(t_{hv}) \to \dot{f}(t)$ at $t = t_0$. □

Now let $\tilde{f}_n^*$ be the left derivative of the least concave majorant of the distribution function of a sample of size $n$ drawn from $\tilde{f}_n$. Computationally, this is easy to do using rejection sampling applied to $\tilde{f}_n$ so that normalization of $\tilde{f}_n$ is not needed. We have the following result:

**Proposition 2.1.** *Under conditions A1–A3 and B1–B4,*

$$n^{1/3}(\tilde{f}_n^*(t_0) - \tilde{f}_n(t_0)) \overset{\mathrm{P}}{\underset{*}{\rightsquigarrow}} c(t_0)\mathbb{C},$$

*where* $*$ *denotes the random component of the smoothed bootstrap. In other words, the proposed smoothed bootstrap is consistent in probability.*

*Sketch of proof.* The proof follows the same general arguments used in the proof of the weak convergence of $n^{1/3}(\hat{f}_n(t) - f(t))$. The main idea is that because $\tilde{f}_n$ satisfies the conclusions of Lemma 2.1, it satisfies assumptions A1 and A2, for all $n$ large enough with probability approaching 1, and both $\tilde{f}_n(t_0) \overset{\mathrm{P}}{\to} f(t_0)$ and $\tilde{f}_n^{(1)}(t_0) \overset{\mathrm{P}}{\to} \dot{f}(t_0)$. The key challenge is to obtain empirical process results for $\tilde{\mathbb{G}}_n \equiv \sqrt{n}(\tilde{\mathbb{P}}_n - P_n)$, where $\tilde{\mathbb{P}}_n$ is the empirical distribution for an i.i.d. sample from $P_n$, where $P_n$ changes with $n$. We need modulus of continuity results which are uniform in $P_n$. Some related uniformity concepts and results are given in Chapter 2.8 of van der Vaart and Wellner [18]. In our case, $P_n$ is the probability measure obtained by integrating $\hat{f}_n$. We omit the remaining details (which are lengthy). $\square$

## 3. $L_1$ Error

Under A1–A2 and the additional assumption

A3′. $\sup_{t \in (0,1)} |\dddot{f}(t)| < \infty$,

Groeneboom et al. [12] proved that

$$n^{1/6} \left\{ n^{1/3} \int_0^1 |\hat{f}_n(t) - f(t)| dt - \mu(f) \right\} \rightsquigarrow N(0, \sigma^2),$$

where

$$\sigma^2 \equiv 8 \int_0^\infty \mathrm{cov}\left(|\xi(0)|, |\xi(x)|\right) dx$$

and, for a differentiable density $g$,

$$\mu(g) \equiv 2\mathrm{E}|\xi(0)| \int_0^1 \left| \frac{1}{2} \dot{g}(t) g(t) \right|^{1/3} dt.$$

In the above, the process $\xi$ is a stationary process constructed from a two-sided Brownian motion $\mathbb{Z}$ as follows:

$$\xi(t) \equiv \operatorname*{arg\,max}_{h \in \mathbb{R}} \left\{ \mathbb{Z}(t + h) - \mathbb{Z}(t) - h^2 \right\}.$$

We will utilize the smoothed bootstrap $\tilde{f}_n^*$ again, for inference in this setting, but we will need a "smoother" kernel and larger bandwidth. In particular, we need the additional assumptions

B5. $\int_{-1}^1 vK(v)dv = 0$, $\int_{-1}^1 \ddot{K}(v)dv = 0$, and $\int_{-1}^1 v\ddot{K}(v)dv = 0$; and

B6. $|(d/(dt))\ddot{K}(t)|$ is uniformly bounded over $(-1, 1)$.

A kernel that satisfies B1–B6 is $K(v) = (15/16)(1 - v^2)^2$.

We have the following lemma. The proof is similar to the proof of Lemma 2.1, and we omit the details.

**Lemma 3.1.** *Assume $h = R_n n^{-\alpha}$, where $0 < R_n + R_n^{-1} = O_P(1)$ and $\alpha \in (1/6, 1/5)$. Under A1–A2, A3′, and B1–B6, we have the following:*

(i) $\sup_{t \in [0,1]} |\tilde{f}_n(t) - f(t)| = O_P(n^{-2\alpha})$;

(ii) $\sup_{t \in (0,1)} |\tilde{f}_n^{(1)}(t) - \dot{f}(t)| = O_P(n^{-\alpha})$; *and*

(iii) *There exists a constant $a < \infty$ such that*

$$\sup_{t \in (0,1)} |\tilde{f}_n^{(2)}(t)| \leq a + o_P(1).$$

*In particular, we have that*

$$\int_0^1 \left| \frac{1}{2} \tilde{f}_n^{(1)}(t) \tilde{f}_n(t) \right|^{1/3} dt - \int_0^1 \left| \frac{1}{2} \dot{f}(t) f(t) \right|^{1/3} dt = o_P(n^{-1/6}).$$

Here is the proposed procedure. Let

- $\tilde{f}_n^*$ be a bootstrapped Grenander estimator based on a i.i.d. sample $X_1^*, \ldots, X_n^*$ drawn from $\tilde{f}_n$ and
- $\tilde{f}_{n,m}^{**}$ be an additional bootstrapped Grenander estimator based on an i.i.d. sample $X_1^{**}, \ldots, X_m^{**}$ also drawn from $\tilde{f}_n$ but independent of the first sample, where $m$ is much larger than $n$, i.e., we require $m/n \to \infty$.

Compute

$$\hat{\mu}_{n,m} = m^{1/3} \int_0^1 \left| \tilde{f}_{n,m}^{**}(t) - \tilde{f}_n \right| dt$$

based on one large bootstrap realization. Finally, estimate the $1 - \alpha$ upper critical value, which estimate we denote $\hat{C}_\alpha$, of the bootstrapped distribution of

$$n^{1/6} \left\{ n^{1/3} \int_0^1 \left| \tilde{f}_n^*(t) - \tilde{f}_n(t) \right| dt - \hat{\mu}_{n,m} \right\},$$

based on repeating the bootstrap $\tilde{f}_n^*$. Note that the large bootstrap (a "supersample bootstrap") is only computed once in the process of obtaining this critical value.

We have the following proposition:

**Proposition 3.1.** *Assume the conditions of Lemma 3.1 and that $\hat{C}_\alpha$ is obtained through using the above procedure. Then the set of densities*

$$\left\{ g : \int_0^1 \left| \hat{f}_n(t) - g(t) \right| dt \leq n^{-1/3} \hat{\mu}_{n,m} + n^{-1/2} \hat{C}_\alpha \right\}$$

*has asymptotically $(1 - \alpha)$ coverage, provided $m/n \to \infty$.*

*Sketch of proof.* The proof is very lengthy and we omit the details. As with the proof of Proposition 2.1, the basic argument is that $\tilde{f}_n$ shares the required properties of $f$ for all sufficiently large $n$ with probability approaching 1, as a consequence of Lemma 3.1. Under these circumstances, the arguments in the proof given in Groeneboom et al. [12] can be carried over from $f$ to $\tilde{f}_n$. Accordingly, we first obtain that

$$n^{1/6} \left\{ n^{1/3} \int_0^1 |\tilde{f}_n^*(t) - \tilde{f}_n(t)| dt - \mu(\tilde{f}_n) \right\} \overset{\mathrm{P}}{\underset{*}{\rightsquigarrow}} N(0, \sigma^2)$$

and

$$m^{1/6} \left\{ m^{1/3} \int_0^1 |\tilde{f}_{n,m}^{**}(t) - \tilde{f}_n(t)| dt - \mu(\tilde{f}_n) \right\} \overset{\mathrm{P}}{\underset{*}{\rightsquigarrow}} N(0, \sigma^2).$$

Thus $\hat{\mu}_{n,m} - \mu(\tilde{f}_n) = O_P(m^{-1/6}) = o_P(n^{-1/6})$, conditionally on the data, and hence both

$$(3.1) \qquad n^{1/6} \left\{ n^{1/3} \int_0^1 |\tilde{f}_n^*(t) - \tilde{f}_n(t)| dt - \hat{\mu}_{n,m} \right\} \quad \overset{\mathrm{P}}{\underset{*}{\rightsquigarrow}} \quad N(0, \sigma^2)$$

and

$$n^{1/6} \left\{ n^{1/3} \int_0^1 |\hat{f}_n(t) - f(t)| dt - \hat{\mu}_{n,m} \right\}$$
$$= \quad n^{1/6} \left\{ n^{1/3} \int_0^1 |\hat{f}_n(t) - f(t)| dt - \mu(f) \right\} + o_P(1)$$
$$\rightsquigarrow \quad N(0, \sigma^2),$$

since $\mu(\tilde{f}_n) = \mu(f) + o_P(n^{-1/6})$ by lemma 3.1 and the restriction that $\alpha > 1/6$. Combining this with (3.1), we obtain the desired conclusion. $\square$

## 4. Discussion and open questions

We note that Abrevaya and Huang [1], in their Theorem 3, provide a general result on the unconditional limiting distribution of the bootstrap for argmax estimators that implies (2.5) and includes many other monotone function settings. Their result, in combination with our arguments and our Theorem 2.2, could thus probably be used to deduce bootstrap inconsistency for monotone function estimators in general.

Note also that under the conditions of Lemma 3.1,

$$\sup_{t \in [0,1]} \left| \tilde{f}_n(t) - f(t) \right| = o_P(n^{-1/3}).$$

This means that the assumptions on the smoothness of $f$ used in Groeneboom et al. [12] are so strong that we can construct a kernel density estimator that uniformly converges faster than the Grenander estimator. Thus, with these assumptions, the Grenander estimator is clearly not optimal. This raises the important question about whether the assumptions in Groeneboom et al. [12] can be relaxed to the point that there are no kernel density estimators superior to the Grenander estimator. Alternatively, is it possible to show that the assumptions cannot be relaxed? If this is the case, then the Grenander is generally not optimal for $L_1$ confidence band construction.

Perhaps a more pressing open problem is constructing valid uniform confidence bands for the Grenander estimator. It appears as if establishing the uniform rate, which seems to be $n^{1/3}(\log n)^{-1/3}$, is not too hard in comparison to establishing distributional convergence. As with the arguments used in Groeneboom et al. [12] for the $L_1$ error, it appears as if the uniform error should converge to some extremum of the process $|\xi(t)|$ defined in section 3 over some increasing interval $[0, \tau_n]$. If this could be done, then the extremal limiting distribution results in Hooghiemstra and Lopuhaä [13] may be applicable, yielding an extreme value distribution in the limit after standardization. Establishing this, however, seems to be very difficult without results for convergence of empirical processes over noncompact index sets. A further question is whether this can be accomplished without imposing assumptions so strong that the primacy of the Grenander is lost (as seems to have happened in the $L_1$ error case). This issue of lost primacy, of course, does not arise in the pointwise error setting (Birgé [4]).

Finally, we note that these results and issues for the Grenander estimator have implications for the survival estimator under current status censoring as well as for monotone function estimation in general because of the similar argmax structures noted previously.

**Acknowledgments.** Thanks to Moulinath Banerjee for many insightful and interesting discussions on the bootstrap and for encouraging me to write this paper. Thanks also to Bodhisattva Sen and Michael Woodroofe for their insightful participation in these bootstrap discussions. Finally, thanks also to an anonymous referee for a very insightful review that led to important improvements in the paper. The paper is based on a contributed talk I gave of the same title at the 2004 Joint Statistical Meetings in Toronto. I also wish to take this opportunity to express appreciation and respect for the life and work of Pranab K. Sen.